\documentclass{amsart}
\usepackage{parskip}
\usepackage{amsmath,amssymb,amsthm,amsfonts,mathrsfs}
\usepackage{fullpage}
\usepackage{xcolor,microtype}
\usepackage[hidelinks]{hyperref}
\allowdisplaybreaks

\usepackage{cite}
\usepackage{colonequals}

\usepackage[capitalize]{cleveref}

\newtheorem{theorem}{Theorem}[section]
\newtheorem{lemma}[theorem]{Lemma}
\newtheorem{corollary}[theorem]{Corollary}
\newtheorem{proposition}[theorem]{Proposition}

\theoremstyle{definition}

\newtheorem{example}[theorem]{Example}
\newtheorem{remark}[theorem]{Remark}
\newcommand{\N}{\mathbb{N}}
\newcommand{\Z}{\mathbb{Z}}

\renewcommand{\P}{\mathbb{P}}
\newcommand{\E}{\mathbb{E}}

\newcommand{\cP}{\mathcal P}

\newcommand{\cL}{\mathcal{L}}
\newcommand{\cA}{\mathcal{A}}
\newcommand{\cB}{\mathcal{B}}
\newcommand{\cC}{\mathcal{C}}

\newcommand{\set}[1]{\left\{#1\right\}}

\newcommand{\Prob}{\mathbb P}
\newcommand{\1}{\mathbf 1}
\newcommand{\eps}{\varepsilon}
\newcommand{\dens}{\operatorname{dens}}
\newcommand{\lcmop}{\operatorname{lcm}}

\newcommand{\diag}{\operatorname{diag}}

\newcommand{\abs}[1]{\lvert #1\rvert}
\newcommand{\cN}{\mathcal N}

\newcommand{\cR}{\mathcal R}

\author[L. Duan]{Lian Duan}
\address{Institute of Mathematical Sciences, ShanghaiTech University, No. 393 Middle Huaxia Road, Pudong New District, Shanghai 201210, China}
\email{duanlian@shanghaitech.edu.cn}

\author[Z. Ge]{Zhenchao Ge}
\address{Department of Pure Mathematics, University of Waterloo, 200 University Avenue West, Waterloo, ON, N2L 3G1, Canada}
\email{z4ge@uwaterloo.ca}

\author[B. Wang]{Biao Wang}
\address{School of Mathematics and Statistics, Yunnan University, Kunming, Yunnan 650500, China}
\email{bwang@ynu.edu.cn}

\author[C. Wang]{Chunlin Wang}
\address{School of Mathematical Sciences, Sichuan Normal University, Chengdu 610064, China}
\email{c-l.wang@outlook.com}

\author[B. Xue]{Boqing Xue}
\address{Institute of Mathematical Sciences, ShanghaiTech University, No. 393 Middle Huaxia Road, Pudong New District, Shanghai 201210, China}
\email{xuebq@shanghaitech.edu.cn}

\date{}

\makeatletter
\@namedef{subjclassname@2020}{\textup{}2020 Mathematics Subject Classification}
\makeatother

\title{On covering systems over finite nilpotent groups and lattices}
\subjclass[2020]{11B25, 11A07}
\keywords{Covering systems, minimum modulus problem, nilpotent groups}

\theoremstyle{plain}

\begin{document}

\begin{abstract}
We prove a quantitative non-covering theorem for finite nilpotent groups. As a consequence, every finite covering of $\Z^d$ by affine full-rank sublattices with pairwise distinct indices contains a member of bounded index, uniformly in $d$. We also establish a summability criterion allowing repeated indices,
which extends to finite groups containing a normal nilpotent subgroup of bounded index.
\end{abstract}

\maketitle

\section{Introduction and statement of results}\label{sec_intro}

In 1950, Erd\H{o}s \cite{Erdos1950} introduced a \textit{covering system} of the integers $\Z$, namely a finite collection of residue classes whose union is $\Z$. He asked whether the minimum modulus of a covering system with distinct moduli can be arbitrarily large. This question became known as the \textit{minimum modulus problem}. In 2015, building on the work of Filaseta, Ford, Konyagin, Pomerance and Yu \cite{FFSPY2007}, Hough \cite{Hough2015} resolved the problem by proving that the minimum modulus is at most $10^{16}$. In 2022, Balister, Bollob\'as, Morris, Sahasrabudhe and Tiba \cite{BBMST2022} reduced this bound to $616{,}000$ and established a quantitative lower bound for the density of the uncovered set when all moduli are sufficiently large. Other developments include results on divisibility restrictions on the moduli \cite{HoughNielsen2019} and on the structure and enumeration of minimal covering systems \cite{BBMST2024}.

Covering systems have also been studied in more general algebraic settings. Covers of groups by cosets of subgroups were investigated by Neumann \cite{Neumann1954} and Sun \cite{Sun2006}. Analogues of the minimum modulus theorem have been proved for rings of integers of number fields by Li, Wang and Yi \cite{LiWangYi2023pre}, and for polynomial rings over finite fields by Li, Wang, Wang and Yi \cite{LiWangWangYi2023pre}. The latter authors also obtained nonexistence results for covering systems in global function fields \cite{LiWangWangYi2025}.

Our first result extends the quantitative uncovered-set estimate to finite nilpotent groups. For any $\eps>0$, we define the multiplicative function $\mu_\eps$ by $\mu_\eps(1)=1$ and
\[
\mu_\eps(p^a)=
\begin{cases}
1,&\text{ if } p=2, \\
1+p^{-1}(\log p)^3 \exp\big((\log\log p)^\eps-1\big), & \text{ if } p\ge3
\end{cases}
\]
for every prime $p$ and every integer $a\geq1$.

\begin{theorem}\label{thm:finite-nilpotent}
For every $\eps>0$, there exists an integer $M=M(\eps){\geq2}$ with the following property. Let $G$ be a finite nilpotent group, and let $C_i=a_iH_i$ $(1\leq i\leq k)$ be left cosets of subgroups $H_i\leq G$. Suppose that the indices $D_i=[G:H_i]$ are pairwise distinct and satisfy $D_i\geq M$. Then
\[
 \frac{\abs{G\setminus\bigcup_{i=1}^kC_i}}{\abs G}
 \geq
 \frac{1}{2}\exp\left(
 -4\sum_{i=1}^k\frac{\mu_\eps(D_i)}{D_i}
 \right).
\]
\end{theorem}

We next apply Theorem~\ref{thm:finite-nilpotent} to affine lattices. Let $M_d(\Z)$ denote the set of $d\times d$ integral matrices. An affine full-rank sublattice of $\Z^d$ has the form
\[
L=t+\Lambda=t+A\Z^d,
\qquad t\in\Z^d,\quad A\in M_d(\Z),\quad \det A\neq0.
\]
Its index is $D=[\Z^d:\Lambda]=\abs{\det A}$. A finite family of such lattices is a covering system of $\Z^d$ if its union is $\Z^d$. Borodzik, Nguyen and Robins \cite{BNR} studied disjoint coverings by translated sublattices, which they call lattice tilings. Cremona and Koymans \cite{CK} studied homogeneous
coverings of $\Z^2$ by finite-index sublattices and systematic refinement constructions. Here we allow both translates and overlaps, and impose conditions on the index multiplicities.

For $S\subseteq\Z^d$, its \textit{asymptotic density} is defined by
\[
\dens(S)=\lim_{N\to\infty}
\frac{|S\cap[-N,N]^d|}{(2N+1)^d},
\]
provided that the limit exists. Given $L_i=t_i+A_i\Z^d$ $(1\leq i\leq k)$, put $D_i=\abs{\det A_i}$ and $Q=\lcmop(D_1,\ldots,D_k)$. Since $Q\Z^d\subseteq D_i\Z^d\subseteq A_i\Z^d$, each $L_i$ is the inverse image of a coset of index $D_i$ in $G=(\Z/Q\Z)^d$. The uncovered set is periodic modulo $Q$, so its asymptotic density exists and equals the uncovered proportion in $G$. As $G$ is abelian, and hence nilpotent, Theorem~ \ref{thm:finite-nilpotent} gives the following corollary.

 \begin{corollary} \label{thm:lattice-main}
Fix $\eps>0$. The constant $M=M(\eps)$ in Theorem~\ref{thm:finite-nilpotent} has the following property for every integer $d\geq1$. Let $L_i=t_i+A_i\Z^d$ $(1\leq i\leq k)$ be affine full-rank sublattices of $\Z^d$. Suppose that the indices $D_i=\abs{\det A_i}$ are pairwise distinct and satisfy $D_i\geq M$.
Then the uncovered set $R=\Z^d\setminus\bigcup_{i=1}^kL_i$ satisfies
\[
 \dens(R)
 \geq
 \frac{1}{2}\exp\left(
 -4\sum_{i=1}^k\frac{\mu_\eps(D_i)}{D_i}
 \right) {>0}.
 \]
 In particular, the family does not cover $\Z^d$.
\end{corollary}

For $d=1$ {and $0<\eps\leq1$}, Corollary ~\ref{thm:lattice-main} gives a refinement of \cite[Theorem~1.1]{BBMST2022}. For $p\geq3$, concavity gives
\[
(\log\log p)^\eps-1
\leq \eps(\log\log p-1)
\leq \eps\log\log p.
\]
Hence, our prime-power weight is no larger than $1+p^{-1}(\log p)^{3+\eps}$. Taking $\eps=1/2$ gives the following consequence.

\begin{corollary}\label{cor:min-index}
There is an absolute constant $M_0$ such that, for every $d\geq1$, every finite covering of $\Z^d$ by affine full-rank sublattices with pairwise distinct indices contains a member of index less than $M_0$.
\end{corollary}

Some restriction on the index multiplicities is necessary: the $D$ cosets of any index-$D$ sublattice partition $\Z^d$. However, pairwise distinctness can be relaxed. Klein, Koukoulopoulos and Lemieux \cite[Theorem~3]{KKL2022} proved that the minimum modulus of an integer covering system is bounded in terms of its multiplicity. We now consider finite nilpotent groups and affine lattices under a more flexible multiplicity condition.

We allow repetitions of both subgroups and cosets in an indexed family. For a finite indexed family $\cC=(a_iH_i)_{i=1}^k$ of left cosets in a finite group $G$, write
\[
n_{\cC}(D)
=
\bigl|\{i\in\{1,\ldots,k\}:[G:H_i]=D\}\bigr|
\]
for the number of members of index $D$, counted with multiplicity. For $D\geq2$, let $P^+(D)$ denote its largest prime factor, $\tau(D)$ its number of positive divisors, and $\omega(D)$ its number of distinct prime factors. We use the conventions $P^+(1)=1$, $\tau(1)=1$ and $\omega(1)=0$.
Let
\[
 m:\{2,3,4,\ldots\}\longrightarrow\N\cup\{0\}
\]
be a function bounding these multiplicities. For each prime $p$, define the following quantities, which control the first and second moments in our argument:
 \begin{align*}
 \mathsf{A}_p(m)
 &:={}
 \sum_{\substack{D\geq2\\P^+(D)=p}}
 \frac{m(D)\tau(D)}{D},
 \\
 \mathsf{B}_p(m)
 &:={}
 \frac{\log^2(2p)}{p}
 \sum_{\substack{D\geq2\\P^+(D)=p}}
 \frac{m(D)^2\tau(D)^2}{D}.
\end{align*}

\begin{theorem}\label{thm:main}
Let $m:\{2,3,4,\ldots\}\longrightarrow\N\cup\{0\}$ satisfy
 \begin{equation}
\label{eq:min-condition}
 \sum_{p} \min\bigl\{\mathsf{A}_p(m),\mathsf{B}_p(m)\bigr\}<\infty.
\end{equation}
Then there exists an integer $M=M(m)\geq2$ with the following property. Let $G$ be a finite nilpotent group, and let $\cC={\{C_i=a_iH_i:1\leq} {i\leq k\}}$ be a finite family of left cosets of subgroups $H_i\leq G$. Suppose that
\[
D_i {=}[{G:H_i}] \geq M {\quad} (1\leq i\leq k),
\qquad n_{\cC}(D)\leq m(D) {\quad} (D\geq2).
\]
Then
 \[
{\frac{|G\setminus\bigcup_{i=1}^kC_i|}{|G|}
\geq 2^{-\omega(Q)-1}>0,
}\]
where $Q=\lcmop(D_1,\ldots,D_k)$. In particular, $\cC$ does not cover $G$.
\end{theorem}

The same finite-quotient argument gives a lattice version. For a family $\cL=\{L_i=t_i+\Lambda_i:1\leq i\leq k\}$, write
\[
n_{\cL}(D)=\lvert\{i:\, [{\Z^d:\Lambda_i}]=D\}\rvert.
\]

 \begin{corollary} \label{thm:main_lattices}
Let $m:\{2,3,4,\ldots\}\longrightarrow\N\cup\{0\}$ satisfy \eqref{eq:min-condition}. The constant $M=M(m)$ in Theorem~\ref{thm:main} has the following property for every integer $d\geq1$. Let $\cL=\{L_i=t_i+\Lambda_i:1\leq i\leq k\}$ be a finite family of affine full-rank sublattices of $\Z^d$. Suppose that 
\[
D_i =[\Z^d:\Lambda_i]\geq M {\quad} (1\leq i\leq k),
\qquad n_{\cL}(D)\leq m(D) {\quad} (D\geq2).
\]
Then the uncovered set $R=\Z^d\setminus\bigcup_{i=1}^kL_i$ satisfies
\[
\dens(R)\geq 2^{-\omega(Q)-1}>0,
\]
where $Q=\lcmop(D_1,\ldots,D_k)$. In particular, $\cL$ does not cover $\Z^d$.
\end{corollary}

\begin{remark}\label{remark1.6}
Condition \eqref{eq:min-condition} holds if and only if the set of primes can be partitioned as $\cP=\cP_1\sqcup\cP_2$ so that
\begin{equation}
\label{eq:first-series}
 \sum_{\substack{D\geq2\\P^+(D)\in\cP_1}}
 \frac{m(D)\tau(D)}{D}<\infty
\end{equation}
and
\begin{equation}
\label{eq:second-series}
 \sum_{\substack{D\geq2\\P^+(D)\in\cP_2}}
 \frac{m(D)^2\tau(D)^2
 \log^2\!\bigl(2P^+(D)\bigr)}{D\,P^+(D)}<\infty.
\end{equation}
\end{remark}

The article is organized as follows. Section~\ref{sec_examp} contains examples illustrating our results. In Section~\ref{sec_distortion_method}, we adapt the distortion method to finite nilpotent groups. We then apply this method to prove Theorems~\ref{thm:finite-nilpotent} and~\ref{thm:main} in Sections~\ref{sec_proof_1} and~\ref{sec_proof_2}, respectively. Finally, in Section~\ref{sec:bounded-index-extension}, we extend Theorem~\ref{thm:main} to finite groups containing a normal nilpotent subgroup of bounded index.

\subsection*{Acknowledgment}

\subsection*{Use of AI}
In this paper, AI tools were used to fill in proof details for results claimed by the authors, who also organized the proof strategy. The authors then verified and simplified the proofs and take full responsibility for the human verification and the content.

\section{ {Examples} }\label{sec_examp}

 We write $e_1,\ldots,e_d$ for the standard basis of $\Z^d$.

 \subsection{{Coverings induced by integer covering systems}}

\begin{example}\label{prop:pullback}
Let $c=(c_1,\ldots,c_d)\in\Z^d$ be primitive, meaning that
$\gcd(c_1,\ldots,c_d)=1$. For each $m\geq1$, define
\[
 \Lambda(c,m)
 =\{x\in\Z^d:c\cdot x\equiv0\pmod m\}.
\]
The homomorphism $\Z^d\longrightarrow\Z/m\Z$, given by $x\longmapsto c\cdot x\pmod m$, is surjective and has kernel $\Lambda(c,m)$. Hence, $[\Z^d:\Lambda(c,m)]=m$.

Suppose that the residue classes $a_i\pmod{m_i}$,
$1\leq i\leq k$, cover $\Z$. Then
\[
 L_i=\{x\in\Z^d:c\cdot x\equiv a_i\pmod{m_i}\},
 \qquad 1\leq i\leq k,
\]
are affine full-rank lattices of indices $m_i$ and cover $\Z^d$.
Indeed, for every $x\in\Z^d$, the integer $c\cdot x$ belongs
to at least one of the given residue classes.
 \end{example}

\begin{example}\label{ex:linear-form}
Take $c=(1,\ldots,1)$. For $d=1$, take $A_m=(m)$. For $d\geq2$, define
\[
 A_m=
 \begin{pmatrix}
 m&-1&-1&\cdots&-1\\
 0& 1& 0&\cdots& 0\\
 0& 0& 1&\ddots& \vdots\\
 \vdots&\vdots&\ddots&\ddots&0\\
 0&0&\cdots&0&1
 \end{pmatrix}.
\]
Then $\det A_m=m$ and
\[
 a e_1+A_m\Z^d
 =\set{x\in\Z^d:x_1+\cdots+x_d\equiv a\pmod m}.
\]
The residue classes
 \[
 0\pmod2,
 \quad0\pmod3,
 \quad1\pmod4,
 \quad5\pmod6,
 \quad7\pmod{12}
\]
 cover $\Z$. Hence,
 \[
 A_2\Z^d,
 \quad A_3\Z^d,
 \quad e_1+A_4\Z^d,
 \quad5e_1+A_6\Z^d,
 \quad7e_1+A_{12}\Z^d
\]
cover $\Z^d$. Their indices are the pairwise distinct integers
$2,3,4,6,12$.
\end{example}

\subsection{A two-dimensional covering system}

\begin{example}\label{prop:photo}
Consider the five affine lattices $L_i=t_i+A_i\Z^2$ given by
\[
\begin{array}{c|c|c|c}
i&A_i&t_i&\abs{\det A_i}\\ \hline
1&
\begin{pmatrix}1&0\\0&2\end{pmatrix}&
\binom{0}{1}&2\\[5pt]
2&
\begin{pmatrix}3&2\\0&1\end{pmatrix}&
\binom{0}{0}&3\\[5pt]
3&
\begin{pmatrix}3&1\\0&2\end{pmatrix}&
\binom{2}{0}&6\\[5pt]
4&
\begin{pmatrix}2&1\\0&2\end{pmatrix}&
\binom{-1}{0}&4\\[5pt]
5&
\begin{pmatrix}6&1\\0&2\end{pmatrix}&
\binom{4}{0}&12.
\end{array}
\]
Their indices are pairwise distinct, and they cover $\Z^2$. Indeed, the first lattice contains every point $(x,y)$ with $y$ odd.
On a line $y=2n$, the second and third lattices contain,
respectively, the points satisfying
\[
 x-n\equiv0\pmod3
 \qquad\text{and}\qquad
 x-n\equiv2\pmod3.
\]
The remaining points satisfy $x-n\equiv1\pmod3$.
The fourth lattice contains those for which $x-n$ is odd.
Any point still remaining therefore satisfies
\[
 x-n\equiv1\pmod3,
 \qquad x-n\equiv0\pmod2.
\]
By the Chinese remainder theorem, $x-n\equiv4\pmod6$,
so the point belongs to the fifth lattice.

This construction extends to every dimension $d>2$. Define $\widehat A_i=A_i\oplus I_{d-2}$ and $\widehat t_i=(t_i,0,\ldots,0)$. Then
\[
 \widehat t_i+\widehat A_i\Z^d
 =(t_i+A_i\Z^2)\times\Z^{d-2}.
\]
These five affine lattices cover $\Z^d$, and their indices
remain $2,3,6,4,12$, respectively.
\end{example}

\subsection{A \texorpdfstring{$d$}{d}-dimensional covering system}

\begin{example}\label{prop:two-direction}
For $d\geq2$, consider the six affine lattices
\begin{align*}
 L_0&=e_1+\diag(2,1,\ldots,1)\Z^d,\\
 L_2&=\diag(2,2,1,\ldots,1)\Z^d,\\
 L_3&=\diag(2,3,1,\ldots,1)\Z^d,\\
 L_4&=e_2+\diag(2,4,1,\ldots,1)\Z^d,\\
 L_6&=5e_2+\diag(2,6,1,\ldots,1)\Z^d,\\
 L_{12}&=7e_2+\diag(2,12,1,\ldots,1)\Z^d.
\end{align*}
The determinants of the diagonal matrices show that their
indices are the pairwise distinct integers $2,4,6,8,12,24$.

These lattices cover $\Z^d$. If $x_1$ is odd, then
$x=(x_1,\ldots,x_d)\in L_0$. If $x_1$ is even, apply the
five-class covering system from Example~\ref{ex:linear-form}
 to $x_2$ . The corresponding class places $x$ in one of
$L_2,L_3,L_4,L_6,L_{12}$.
\end{example}

\subsection{An example involving both moments}
\label{subsec:mixed-moments}

The following example shows that condition \eqref{eq:min-condition}
can hold even when both full series
$\sum_p\mathsf{A}_p(m)$ and $\sum_p\mathsf{B}_p(m)$ diverge.

\begin{example}\label{ex:mixed-moments}
Define
\[
m(D)=
\begin{cases}
2^j,&D=2^{2j}\text{ for some integer }j\geq2,\\
1,&D\text{ is an odd prime},\\
0,&\text{otherwise}.
\end{cases}
\]
Since $\tau(2^{2j})=2j+1$, we have
\[
\mathsf{A}_2(m)
=\sum_{j\geq2}\frac{2j+1}{2^j}<\infty,
\qquad
\mathsf{B}_2(m)
=\frac{\log^2 4}{2}\sum_{j\geq2}(2j+1)^2=\infty.
\]
For every odd prime $p$, we have $\mathsf{A}_p(m)=2/p$ and $\mathsf{B}_p(m)=4\log^2(2p)/p^2$. Consequently,
\[
\sum_p\mathsf{A}_p(m)=\infty,
\qquad
\sum_p\mathsf{B}_p(m)=\infty,
\]
whereas
\[
\sum_p\min\{\mathsf{A}_p(m),\mathsf{B}_p(m)\}
\leq
\mathsf{A}_2(m)
+\sum_{\substack{p\geq3\\p\ {\rm prime}}}
\frac{4\log^2(2p)}{p^2}
<\infty.
\]
Thus, neither full summability condition holds, but the combined
condition \eqref{eq:min-condition} does. In the partition formulation
\eqref{eq:first-series}--\eqref{eq:second-series}, we take
\[
\cP_1=\{2\},
\qquad
\cP_2=\{p:p\text{ is an odd prime}\}.
\]
 Theorem ~\ref{thm:main} and {Corollary~} \ref{thm:main_lattices}
therefore apply to this multiplicity function.

For a concrete lattice construction, fix $d\geq1$, integers
 $N\geq {J}\geq2$, and a finite set $S$ of odd primes.
For ${J}\leq j\leq N$, let
\[
\cA_j
=\{a\in\Z:0\leq a<2^{2j},\
  2^{j-1}\mid a,\ 2^j\nmid a\}.
\]
There are exactly $2^j$ elements in $\cA_j$. Define
\[
L_{j,a}
=\{x=(x_1,\ldots,x_d)\in\Z^d:
 x_1\equiv a\pmod{2^{2j}}\},
\qquad a\in\cA_j.
\]
For each $p\in S$, choose an arbitrary integer $b_p$ and put
\[
L_p
=\{x=(x_1,\ldots,x_d)\in\Z^d:
 x_1\equiv b_p\pmod p\}.
\]
The family
\[
\cL
=\{L_{j,a}: {J} \leq j\leq N,\ a\in\cA_j\}
 \cup\{L_p:p\in S\}
\]
has $2^j$ members of index $2^{2j}$ for each ${J}\leq j\leq N$
 and one member of index $p$ for each $p\in S$.
Hence, its index multiplicities are bounded by $m$.

For each $j$, the union of the lattices $L_{j,a}$ over
$a\in\cA_j$ is
\[
\{x\in\Z^d:2^{j-1}\mid x_1,\ 2^j\nmid x_1\},
\]
which has density $2^{-j}$. These unions are pairwise disjoint
as $j$ varies. Since their union is periodic modulo $2^{2N}$,
the Chinese remainder theorem gives
\[
 \begin{aligned}
\dens\left(\Z^d\setminus\bigcup_{L\in\cL}L\right)
=\left(1-\sum_{j=J}^{N}2^{-j}\right)
 \prod_{p\in S}\left(1-\frac{1}{p}\right) =\left(1-2^{1-J}+2^{-N}\right)
 \prod_{p\in S}\left(1-\frac{1}{p}\right)>0.
\end{aligned}
\]
All indices can be made arbitrarily large by choosing ${J}$
 and the primes in $S$ sufficiently large.
\end{example}

\subsection{An example in a nonabelian nilpotent group}

\begin{example}\label{ex:dihedral-covering}
Let
\[
D_8=\langle {\rho}, {\sigma}: {\rho} ^4= {\sigma} ^2=1,\ {\sigma\rho\sigma} = {\rho} ^{-1}\rangle
\]
be the dihedral group of order $8$, and put $T=\langle{\sigma}\rangle$.
The group $D_8$ is nilpotent, while $T$ is a nonnormal subgroup
of index $4$.

Fix $\eps>0$, and choose integers $N\geq {n_0}\geq2$ such that
 $4{n_0}\geq M(\eps)$, where $M(\eps)$ is the constant in
Theorem~\ref{thm:finite-nilpotent}. Put
\[
Q_0=\lcmop( {n_0}, {n_0} +1,\ldots,N),
\qquad
G=D_8\times(\Z/Q_0\Z).
\]
Then $G$ is a finite nonabelian nilpotent group.
For each ${n_0}\leq n\leq N$, define
\[
H_n=T\times n(\Z/Q_0\Z),
\]
and choose an arbitrary left coset $C_n=g_nH_n$ in $G$.
Since $n\mid Q_0$, we have
\[
[G:H_n]=[D_8:T]\,
[\Z/Q_0\Z:n(\Z/Q_0\Z)]=4n.
\]
Thus, the indices are pairwise distinct and at least $M(\eps)$,
although none of the subgroups $H_n$ is normal.

Since $\mu_\eps(2^a)=1$ for every $a\geq1$, we have
$\mu_\eps(4n)=\mu_\eps(n)$. Theorem~\ref{thm:finite-nilpotent}
therefore gives
\[
 {\frac{\left|G\setminus\bigcup_{n=n_0}^{N}C_n\right|}{|G|}}
 \geq
\frac{1}{2}\exp\left(
-\sum _{{{n=n_0}}} ^{N}\frac{\mu_\eps(n)}{n}
\right)>0.
\]
In particular, these cosets cannot cover $G$, regardless of
the choices of their representatives.

On the other hand,
\[
\sum _{{{n=n_0}}} ^{N}\frac{|C_n|}{|G|}
=\frac{1}{4}\sum _{{{n=n_0}}} ^{N}\frac{1}{n},
\]
which can exceed $1$ when $N$ is sufficiently large.
Thus, the theorem gives a positive lower bound even when
the elementary union bound does not.
\end{example}

\section{The distortion method}\label{sec_distortion_method}

Inspired by Hough's work \cite{Hough2015}, Balister et al. \cite{BBMST2022} developed the distortion method to study covering systems over $\Z$. We adapt this method to finite nilpotent groups.
Let $G$ be a finite nilpotent group. For each prime $p\mid |G|$, let $G_p$ denote the unique Sylow $p$-subgroup of $G$. Since a finite group is nilpotent if and only if it is the direct product of its Sylow subgroups, i.e., $G=\prod_p G_p$ (see \cite[Proposition~II.7.5]{Hungerford1980}), we have the following analogue of the Chinese remainder theorem.

\begin{lemma}\label{lem:primary}
Every subgroup $H\leq G$ decomposes as $H=\prod_p(H\cap G_p)$. If $C=aH$ is a coset of $H$ in $G$, then $C=\prod_p C_p$ for suitable cosets $C_p$ of $H\cap G_p$, and
\[
 [G:H]=\prod_p[G_p:H\cap G_p].
\]
\end{lemma}

\begin{proof}
For every $h\in H$, its $p$-component is a power of $h$
and therefore belongs to $H$. Thus,
$H=\prod_p(H\cap G_p)$.
Writing $a=\prod_p a_p$, we obtain
$aH=\prod_p a_p(H\cap G_p)$.
The index identity follows by taking cardinalities.
This completes the proof.
\end{proof}

Fix a finite family of cosets $C_i=a_iH_i\subseteq G$, where $H_i<G$, and put $D_i=[G:H_i]\geq2$. Let
\[
 Q=\lcmop(D_1,\ldots,D_k)
 =\prod_{j=1}^s p_j^{\gamma_j},
 \qquad (p_1<\cdots<p_s).
\]
For each prime $p\mid |G|$ with $p\nmid Q$, the index $[G_p:H_i\cap G_p]$ is a power of $p$ dividing $D_i$, so it equals
$1$ for every $i$. Thus, the normal subgroup
\[
K=\prod_{\substack{p\mid |G|\\p\nmid Q}}G_p
\]
is contained in every $H_i$. Passing to $G/{K}$ preserves the indices
and the uncovered proportion. We may therefore assume that $G=G_1\times\cdots\times G_s$,
where $G_j=G_{p_j}$ is the Sylow $p_j$-subgroup. Each $p_j$ divides at least one of the indices $D_i$.

Set
\[
 G_{<j}=\prod_{\ell<j}G_\ell, \qquad G_{\leq j}=\prod_{\ell\leq j}G_\ell, \qquad G_{>j}=\prod_{\ell>j}G_\ell.
\]
We also set 
\begin{equation}
\label{def:Nj}
\cN_j:=\big\{i: P^+(D_i)=p_j\big\}, \qquad \cB_j:=\bigcup_{i\in\cN_j}C_i.
\end{equation}
The sets $\cN_j$ partition $\set{1,2,\dots,k}$. For $i\in\cN_j$, $C_i$ is the preimage of a coset under the projection $G\to G_{\leq j}$. We identify $C_i$ and $\cB_j$ with their images in $G_{\leq j}$ when no confusion can arise.

\begin{lemma}\label{lem:fibre}
Fix $j$ and let $i\in\cN_j$. Decompose
\[
 D_i= g_i p_{j} ^{r_i},
 \qquad \text{with }p_j\nmid g_i, \text{ and }r_i\geq1.
\]
Under the product decomposition
$G_{\leq j}=G_{<j}\times G_j$, $C_i$ (in $G_{\leq j}$) has the form $E_i\times F_i$, where $E_i$ is a coset of index $g_i$ in $G_{<j}$ and $F_i$ is a coset of index $p_j^{r_i}$ in $G_j$. Therefore, for every $x\in G_{<j}$,
\[
 \frac{\abs{\set{y\in G_j:(x,y)\in C_i}}}{\abs{G_j}}
 =p_{j}^{-r_i}\1_{E_i}(x),
\]
where $\1_{E_i}$ denotes the indicator function of $E_i$. In particular, every nonempty fibre has relative size $p_j^{-r_i}$.
\end{lemma}
\begin{proof}
By Lemma~\ref{lem:primary} and $P^+(D_i)=p_j$, we have
$C_i=E_i\times F_i\times G_{>j}$, where $E_i$ and $F_i$
are cosets of indices $g_i$ and $p_j^{r_i}$ in
$G_{<j}$ and $G_j$, respectively.
For $x\in G_{<j}$, the fibre is $F_i$ if $x\in E_i$
and is empty otherwise.
Since $|F_i|/|G_j|=p_j^{-r_i}$, the asserted identity follows.
This completes the proof.
\end{proof}

We now adapt the distortion method of \cite{BBMST2022} to finite nilpotent groups. We record the measure construction and the estimates needed below.

Let $\Prob_0$ be the
uniform probability measure on $G$. Choose parameters $0\leq\delta_j\leq\frac{1}{2}$ for $1\leq j\leq s$. For $x\in G_{<j}$, define
\begin{equation}\label{eq:alpha}
 \alpha_j(x)
 =\frac{\abs{\set{y\in G_j:(x,y)\in \cB_j}}}{\abs{G_j}}.
\end{equation}
We use the same notation for a probability measure on $G$
and its marginals. At stage $j$, we define $\Prob_j$ on
$G_{\leq j}$ and extend it uniformly over $G_{>j}$.
In particular, before stage $j$,
\[
\Prob_{j-1}(x,y)=\frac{\Prob_{j-1}(x)}{|G_j|}
\qquad(x\in G_{<j},\ y\in G_j).
\]
 Given $\Prob_{j-1}$, define $\Prob_j$ fibre by fibre by
\begin{equation}\label{eq:distorted-measure}
\Prob_j(x,y)=
\begin{cases}
\displaystyle
\max\left\{0,
\frac{\alpha_j(x)-\delta_j}{\alpha_j(x)(1-\delta_j)}
\right\}\Prob_{j-1}(x,y),
&(x,y)\in \cB_j,\\[1.4em]
\displaystyle
\min\left\{
\frac{1}{1-\alpha_j(x)},
\frac{1}{1-\delta_j}
\right\}\Prob_{j-1}(x,y),
&(x,y)\notin \cB_j.
\end{cases}
\end{equation}
If $\alpha_j(x)=0$, the first case is vacuous; if $\alpha_j(x)=1$, the second case is vacuous. The construction preserves the marginal on $G_{<j}$. Consequently,
\begin{equation}\label{invariant_prob}
	\Prob_j(\cB_j)= \Prob_{j+1}(\cB_j)=\Prob_{j+2}(\cB_j)=\cdots=\Prob_{s}(\cB_j)
\end{equation}
for all $1\leq j\leq s$.

Define
\[
 M_j^{(1)}=\E_{j-1}[\alpha_j],
 \qquad
 M_j^{(2)}=\E_{j-1}[\alpha_j^2]
\]
to be the first and second moments of $\alpha_j$ with respect to $\Prob_{j-1}$. Whenever $\delta_j=0$, we interpret $M_j^{(2)}/(4\delta_j(1-\delta_j))$ as $+\infty$. Define a multiplicative function $\nu$ by
\begin{equation}\label{eq:nu}
 \nu( n)=\prod_{{p_j\mid n}} \frac{1}{1-\delta_j}.
\end{equation}
Hence $\nu(p_j^a)=(1-\delta_j)^{-1}$ for every integer $a\geq1$.

We record two estimates which come from adapting \cite[Lemmas~3.3 and~3.4]{BBMST2022}. The construction gives
\[
\Prob_j(\cB_j)
=\E_{j-1}\!\left[
\frac{\max\{\alpha_j-\delta_j,0\}}{1-\delta_j}
\right].
\]
For $0\leq t\leq1$ and $0\leq\delta\leq1/2$, the quantity $\max\{t-\delta,0\}/(1-\delta)$ is at most $t$. If $\delta>0$,
it is also at most $t^2/(4\delta(1-\delta))$. Therefore,
\begin{equation}\label{lem:p_j-upper-bound}
\Prob_j(\cB_j)\leq
\min\left\{M_j^{(1)},
\frac{M_j^{(2)}}{4\delta_j(1-\delta_j)}\right\}.
\end{equation}

For a coset $C=aH$ of index $h=[G:H]$, we also have
\begin{equation}\label{lem:coset-measure}
\Prob_j(C)\leq \frac{1}{h}\prod_{\substack{\ell\leq j\\p_\ell\mid h}}
\frac{1}{1-\delta_\ell}\leq\frac{\nu(h)}{h}. 
\end{equation}
Indeed, write $C=\prod_\ell C_\ell$ as in Lemma~\ref{lem:primary}. If $p_j\nmid h$, then $C_j=G_j$, so preservation of the marginal on $G_{<j}$ and uniformity on $G_{>j}$ give $\Prob_j(C)=\Prob_{j-1}(C)$. If $p_j\mid h$, the pointwise bound $\Prob_j\leq(1-\delta_j)^{-1}\Prob_{j-1}$ gives the corresponding inequality for $C$. Induction from $\Prob_0(C)=1/h$ proves \eqref{lem:coset-measure}.

The next lemma gives a lower bound for the uncovered proportion.

\begin{lemma}[{\cite[Theorem~3.1]{BBMST2022}}]
\label{thm:general-sieve}
Set
$
 \eta
 =\sum\limits_{j=1}^s
 \min\left\{
 M_j^{(1)},
 \frac{M_j^{(2)}}{4\delta_j(1-\delta_j)}
 \right\}.
$
If $\eta<1$, then
\begin{equation}\label{eq:general-sieve}
 \frac{\abs{G\setminus\bigcup_iC_i}}{\abs G}
 \geq
 (1-\eta)
 \exp\left(
 -\frac{2}{1-\eta}
 \sum_{i=1}^k\frac{\nu(D_i)}{D_i}
 \right).
\end{equation}
\end{lemma}

\begin{proof}
Let $R=G\setminus\bigcup_iC_i$ and $t=\Prob_s(R)$. By
\eqref{invariant_prob}, \eqref{lem:p_j-upper-bound}, and the union bound,
\[
t\geq1-\sum_{j=1}^s\Prob_s(\cB_j)=1-\sum_{j=1}^s\Prob_j(\cB_j)\geq1-\eta>0.
\]
We regard each $\alpha_j$ as a function on $G$ through
projection onto $G_{<j}$, and put
$F=2\sum_{j=1}^s\alpha_j$. 
The multiplier in \eqref{eq:distorted-measure} is at most
$(1-\min\{\alpha_j,\delta_j\})^{-1}$ in either nonvacuous case.
Since $-\log(1-u)\leq2u$ for $0\leq u\leq1/2$, it follows that
\[
\Prob_j(z)\leq e^{2\alpha_j(z)}\Prob_{j-1}(z),
\qquad
\Prob_s(z)\leq e^{F(z)}\Prob_0(z)
\qquad(z\in G).
\]
These inequalities also hold at points of zero mass. Preservation of the earlier marginals gives
\[
\E_s[F]=2\sum_{j=1}^s\E_{j-1}[\alpha_j]
=2\sum_{j=1}^s\Prob_{j-1}(\cB_j)\leq2\sum_{i=1}^k\frac{\nu(D_i)}{D_i},
\]
where the last step uses \eqref{lem:coset-measure}. Consequently, Jensen's inequality and $F\geq0$ imply
\[
\Prob_0(R)\geq\E_s[\1_R e^{-F}]\geq t\exp\bigl(-\E_s[F\mid R]\bigr)\geq t\exp\left(-\frac{2}{t}\sum_{i=1}^k\frac{\nu(D_i)}{D_i}\right).
\]
The function $u\mapsto u e^{-A/u}$ is increasing for $u>0$ when $A\geq0$. Using $t\geq1-\eta$ yields \eqref{eq:general-sieve}.
\end{proof}

\section{Proof of Theorem~\ref{thm:finite-nilpotent}} \label{sec_proof_1}

To prove Theorem~\ref{thm:finite-nilpotent}, we use the following bound for the second moment.

\begin{lemma}[{\cite[Theorem~3.2]{BBMST2022}}]
\label{thm:moments}
Suppose that the indices $D_1,\ldots,D_k$ are pairwise
distinct. For each $j$,
\[
M_j^{(2)} \leq \frac{1}{(p_j-1)^2} \prod_{\ell<j}\left(1+\frac{3p_\ell-1}{(1-\delta_\ell)(p_\ell-1)^2}\right).
\]
\end{lemma}

\begin{proof}
For $i\in\cN_j$, use the notation
$D_i=g_i p_j^{r_i}$ and $C_i=E_i\times F_i\times G_{>j}$
from Lemma~\ref{lem:fibre}. The union bound gives
\[
\alpha_j(x)\leq\sum_{i\in\cN_j}p_j^{-r_i}\1_{E_i}(x).
\]
Consider $i,i'\in\cN_j$ and put $L=\lcmop(g_i,g_{i'})$.
If $E_i\cap E_{i'}$ is nonempty, write it as $bK$, where $K$
is the intersection of their defining subgroups, and put
$h=[G_{<j}:K]$. Then $L\mid h$. By Lemma~\ref{lem:primary}, a
Sylow factor is contained in $K$ exactly when it is contained in
both subgroups. Thus, $h$ and $L$ have the same prime divisors. Hence,
$\nu(h)=\nu(L)$, and \eqref{lem:coset-measure} yields
\[
\Prob_{j-1}(E_i\cap E_{i'})
\leq\frac{\nu(L)}{L}.
\]
This inequality is also valid when the intersection is empty.
Thus,
\[
M_j^{(2)}\leq
\sum_{i,i'\in\cN_j}
\frac{\nu(\lcmop(g_i,g_{i'}))}{p_j^{r_i+r_{i'}}\lcmop(g_i,g_{i'})}.
\]
Distinctness of the $D_i$ means that each pair $(g_i,r_i)$ occurs
at most once. We may therefore extend the sum to all positive
$r,r'$ and all $g,g'$ whose prime divisors lie among
$p_1,\ldots,p_{j-1}$. For a prime $p_\ell$ and $a\geq1$, there
are $2a+1$ pairs of nonnegative exponents with maximum $a$.
Factoring the resulting sum over these primes gives
\[
\begin{aligned}
M_j^{(2)}
&\leq\frac{1}{(p_j-1)^2}
\prod_{\ell<j}\left(
1+\frac{1}{1-\delta_\ell}
\sum_{a=1}^\infty\frac{2a+1}{p_\ell^a}
\right)\\
&=\frac{1}{(p_j-1)^2}
\prod_{\ell<j}\left(
1+\frac{3p_\ell-1}{(1-\delta_\ell)(p_\ell-1)^2}
\right)
\end{aligned}
\]
as desired.
\end{proof}

We also need the following convergence estimate.

\begin{lemma}\label{lem_convergence}
For any $\eps>0$, the following series converges:
\[
\sum_{p\ge3} \frac{1}{p\exp\{(\log\log p)^\eps\}}.
\]
\end{lemma}
\begin{proof}
Let $p_n$ be the $n$th prime. The prime number theorem gives $p_n\gg n\log n$, and $p_n\geq n$. It suffices to
prove the convergence of
\[
\sum_{n\geq3}\frac{1}{n\log n\exp\{(\log\log n)^\eps\}}.
\]
By the integral test and the substitution $y=\log\log t$,
\[
\int_3^\infty\frac{dt}{t\log t\exp\{(\log\log t)^\eps\}}
=\int_{\log\log3}^\infty e^{-y^\eps}\,dy<\infty.
\]
The last integral converges because $e^{-y^\eps}\leq y^{-2}$
for all sufficiently large $y$.
\end{proof}

We follow the proof of \cite[Theorem~1.1]{BBMST2022}, with a different choice of distortion parameters.

\begin{proof}[Proof of Theorem~\ref{thm:finite-nilpotent}]
Fix $\eps>0$ and, for each prime $p\geq3$, put
$
w(p):= (\log p)^3 \exp\!\bigl((\log\log p)^\eps-1\bigr).
$
Let $P_*\geq3$ be a parameter to be chosen later, and set
\begin{equation}\label{eq:delta-choice}
 \delta_j=
 \begin{cases}
 0,&p_j\leq P_*,\\[4pt]
 \displaystyle\frac{w(p_j)}{p_j+w(p_j)},&p_j>P_*.
 \end{cases}
\end{equation}
Since $w(p)=p^{o(1)}$, we may require that $w(p)\leq p$ for every $p>P_*$. Thus, all $\delta_j$ lie in $[0,1/2]$. We use these parameters in \eqref{eq:distorted-measure}.

For $p_j\leq P_*$, we have $\delta_\ell=0$ for every
$\ell<j$, so $\Prob_{j-1}=\Prob_0$.
Therefore, the union bound and distinctness of the indices give
\[
 \sum_{p_j\leq P_*}M_j^{(1)}
 \leq
 \sum_{\substack{i:\ P^+(D_i)\leq P_*}}\frac{1}{D_i} \leq \sum_{\substack{n\geq M\\P^+(n)\leq P_*}}{\frac{1}{n}}
 \leq
 M^{-1/2}
 \prod_{p\leq P_*}(1-p^{-1/2})^{-1}.
\]

Thus, for fixed $P_*$ and sufficiently large $M$, 
\begin{equation}\label{eq:small-quarter}
 \sum_{p_j\leq P_*}M_j^{(1)}\leq\frac{1}{4}.
\end{equation}

Now we consider the case $p_j>P_*$. First applying Lemma~\ref{thm:moments} and using $\log(1+x)\leq x$, we obtain
\[
 M_j^{(2)}
 \leq
 \frac{1}{(p_j-1)^2}
 \exp\left(
 \sum_{\ell<j}
 \frac{3p_\ell-1}{(1-\delta_\ell)(p_\ell-1)^2}
 \right).
\]
By \eqref{eq:delta-choice}, one has $(1-\delta_\ell)^{-1}\leq1+w(p_\ell)/p_\ell$ for every $p_\ell\geq3$.
Consequently, uniformly in $P_*$, for sufficiently large $p_\ell$,
\[
 \frac{3p_\ell-1}{(1-\delta_\ell)(p_\ell-1)^2}
 =\frac{3}{p_\ell}+O_\eps\left(
 \frac{w(p_\ell)}{p_\ell^2}
 \right).
\]
The error terms are summable because $w(p)=p^{o(1)}$. Mertens' estimate gives
\[
 \sum_{p<p_j}\frac{1}{p}=\log\log p_j+O(1).
\]
Consequently, there is a constant $C_\eps$, independent of $P_*$, such that
\begin{equation}\label{eq:log3-bound}
 M_j^{(2)}
 \leq C_{\eps} \frac{(\log p_j)^3}{p_j^2}
 \qquad(p_j>P_*).
\end{equation}
Moreover, since $w(p_j)\leq p_j$ for $p_j>P_*$,
 \[
 \frac{1}{4\delta_j(1-\delta_j)}
 =\frac{(p_j+w(p_j))^2}{4p_jw(p_j)}
 \leq\frac{p_j}{w(p_j)}.
\]
Combining this with \eqref{eq:log3-bound} yields
\[
 \frac{M_j^{(2)}}{4\delta_j(1-\delta_j)}
 \leq
 \frac{e C_\eps}{p_j\exp\{(\log\log p_j)^\eps\}}.
\]
By Lemma~\ref{lem_convergence}, the prime series
\[
\sum_{p\geq 3}\frac{1}{p\exp\{(\log\log p)^\eps\}}
\]
converges. Thus, we may choose $P_*$ sufficiently large so that
\[
 \sum_{p_j>P_*}
 \frac{M_j^{(2)}}{4\delta_j(1-\delta_j)}
 \leq\frac{1}{4}.
\]

Now fix such a $P_*$, and then choose $M=M(\eps)$ so that
\eqref{eq:small-quarter} holds. The parameter $\eta$ in
Lemma~\ref{thm:general-sieve} satisfies $\eta\leq1/2$. Therefore, 
\begin{equation}\label{key_distortion_inequality}
	 \frac{\abs{G\setminus\bigcup_iC_i}}{\abs G}
 \geq
 \frac{1}{2}
 \exp\left(-4\sum_i\frac{\nu(D_i)}{D_i}\right).
\end{equation}

For $p_j>P_*$, the definitions of $\nu$ and $\mu_\eps$, together with \eqref{eq:delta-choice}, give
\[
 \nu(p_j^a)
 =1+\frac{w(p_j)}{p_j}
 =\mu_\eps(p_j^a),
\]
while for $p_j\leq P_*$, one has $\nu(p_j^a)=1\leq\mu_\eps(p_j^a)$.
Both functions are multiplicative, thus $\nu(n)\leq\mu_\eps(n)$ for all
 $n$, and
\[
\sum_i\frac{\nu(D_i)}{D_i} \leq \sum_i\frac{\mu_\eps(D_i)}{D_i}.
\]
Combining this inequality with \eqref{key_distortion_inequality} gives the claimed bound. This finishes the proof. 
\end{proof}

\section{Proof of Theorem~\ref{thm:main}}\label{sec_proof_2}

We use the same notation of Section \ref{sec_distortion_method}
and take 
$\delta_j=\frac{1}{2}$ for $1\leq j\leq s$. It is immediate that inequality \eqref{lem:p_j-upper-bound} implies
\begin{equation}\label{def:M_j^12}
 \Prob_j(\cB_j)
 \leq\min\{M_j^{(1)},M_j^{(2)}\}.
\end{equation}
Assume that $D_i\geq M$ for every $i$, where $M$ will be specified later in the proof.
By \eqref{eq:distorted-measure}, we have $\Prob_j(x,y)\leq2\Prob_{j-1}(x,y)$, and hence for any subset $A\subseteq G$,
\begin{equation}
\label{eq:RN-bound}
\Prob_{s}(A)
\leq 2^s\frac{|A|}{|G|}.
\end{equation}
By \eqref{eq:nu} with $\delta_j=\frac{1}{2}$, we note that $\nu(n)=2^{\omega(n)}$ if every prime factor of $n$ divides $Q$. Thus, using \eqref{lem:coset-measure}, for every coset $C=aH$ in $G$ with index $h=[G:H]$, we have
\begin{equation}\label{lem:coset-mass}
 	 \Prob_j(C)\leq \frac{2^{\omega(h)}}{h}.
 \end{equation}

Fix $j$ and let $\mathcal{N}_j$ be the set defined in \eqref{def:Nj}. For $i\in \mathcal{N}_j$, if we write $D_i=g_i p_j^{r_i}$ with $p_j\nmid g_i$ and $r_i\geq1$, then $P^+(g_i)<p_j$. By Lemma~\ref{lem:fibre}, the image of $C_i$ in $G_{\leq j}=G_{<j}\times G_j$ is $E_i\times F_i$, where $E_i$ is a coset of a subgroup $K_i\leq G_{<j}$
of index $g_i$, and $F_i$ is a coset of index $p_j^{r_i}$ in $G_j$. By \eqref{eq:alpha} and the union bound, we have
\begin{equation}
\label{eq:alpha-upper}
\alpha_j(x)
\leq
\sum_{i\in \cN_j}
\frac{|\{y\in G_j:(x,y)\in C_i\}|}{|G_j|}
=
\sum_{i\in \cN_j}
\1_{E_i}(x)
\,
p_j^{-r_i}.
\end{equation}
Integrating both sides of \eqref{eq:alpha-upper} with respect to $\Prob_{j-1}$
 and using \eqref{lem:coset-mass}, we obtain
\begin{equation}
\label{eq:first-moment-stage}
M_j^{(1)}
\leq
\sum_{i\in \cN_j}
\sum_{x\in E_i}
\Prob_{j-1}(x)
p_j^{-r_i}
\leq
\sum_{i\in \cN_j}
\frac{2^{\omega(g_i)}}{g_i}
p_j^{-r_i}
\leq
\sum_{\substack{D\geq M\\P^+(D)=p_j}}
n_{\cC}(D)\frac{\tau(D)}{D}.\end{equation}
Here, the last inequality follows by grouping the terms with the same index
$D_i=g_ip_j^{r_i}$ and using $2^{\omega(g_i)}\leq\tau(D_i)$.

Next, we consider the second moment. For $i,i'\in\mathcal{N}_j$, when $E_i\cap E_{i'}$ is nonempty, it is a coset of
$K_i\cap K_{i'}$ in $G_{<j}$ with index 
\[
h_{ii'}=[G_{<j}: K_{i}\cap K_{i'}]
=[G_{<j}:K_i][K_i: K_i\cap K_{i'}]
=[G_{<j}:K_{i'}][K_{i'}: K_i\cap K_{i'}].
\]
Consequently, $\lcmop(g_i,g_{i'})\mid h_{ii'}$.
Since $2^{\omega(x)}/x\geq 2^{\omega(y)}/y$ whenever $x\mid y$, \eqref{lem:coset-mass} gives \begin{equation}\label{eq:P(EE)}
 \Prob_{j-1}(E_i\cap E_{i'})
 \leq
\frac{2^{\omega(\lcmop(g_i,g_{i'}))}}{\lcmop(g_i,g_{i'})}.
\end{equation}
When $E_i\cap E_{i'}$ is empty we still have the same bound.
Integrating $\alpha_j(x)^2$ with respect to $\mathbb{P}_{j-1}$ and applying inequalities \eqref{eq:alpha-upper} and \eqref{eq:P(EE)}, we arrive at
\begin{equation}
\label{eq:second-before-group}
M_j^{(2)}
\leq
\sum_{i,i'\in \cN_j}
p_j^{-r_i-r_{i'}}
\Prob_{j-1}(E_i\cap E_{i'})
\leq
\sum_{i,i'\in \cN_j}
p_j^{-r_i-r_{i'}}
\frac{2^{\omega(\lcmop(g_i,g_{i'}))}}{\lcmop(g_i,g_{i'})}.
\end{equation}

To estimate the right-hand side of \eqref{eq:second-before-group},
we establish the following bound for the quadratic form
with kernel $2^{\omega(\mathrm{lcm}(g,h))}/\mathrm{lcm}(g,h)$.

\begin{lemma}
\label{lem:kernel}
There exists an absolute constant $C_0>0$ such that, for every prime $p$
and every finitely supported family of nonnegative numbers
$(x_g)_{P^+(g)<p}$,
\begin{equation}
\label{eq:kernel}
 \sum_{\substack{g,h\geq1\\P^+(g),P^+(h)<p}}
 x_gx_h\,
 \frac{2^{\omega(\lcmop(g,h))}}{\lcmop(g,h)}
 \leq
 C_0\log^2(2p)
 \sum_{P^+(g)<p}
 x_g^2\frac{2^{\omega(g)}\tau(g)}{g}.
\end{equation}
\end{lemma}

\begin{proof}
Using $2x_gx_h\leq x_g^2+x_h^2$, the left-hand side of
\eqref{eq:kernel} is bounded by
\begin{equation}\label{def:R_p}
\sum_{P^+(g)<p}x_g^2\,\cR_p(g),
\quad
\textup{with}
\quad
\cR_p(g):=
\sum_{P^+(h)<p}\frac{2^{\omega(\lcmop(g,h))}}{\lcmop(g,h)}.
\end{equation}
Writing $g=\prod_{q<p}q^{a_q}$ and
summing over the
prime exponents of $h$, we express $\mathcal{R}_p(g)$ as
\[
\cR_p(g) = \prod_{q<p}\Big(\sum_{j\geq 0}\frac{2^{\1_{\max(a_q,j)\geq 1}}}{q^{\max(a_q,j)}}
\Big).
\]
When $a_q=0$ we observe
\[
\sum_{j\geq 0}\frac{2^{\1_{\max(a_q,j)\geq 1}}}{q^{\max(a_q,j)}}
=
1+2\sum_{j\geq 1}q^{-j}
=
1+\frac{2}{q-1},
\]
and when $a_q\geq 1$ we note
\[
\sum_{j\geq 0}\frac{2^{\1_{\max(a_q,j)\geq 1}}}{q^{\max(a_q,j)}}
=
2
\Big(
\sum_{0\leq j\leq a_q}
q^{-a_q}
+
\sum_{j>a_q}
q^{-j}
\Big)
=
\frac{2}{q^{a_q}}\Big(a_q+1+\frac{1}{q-1}\Big).
\]
Combining these identities implies
\[
\cR_p(g)
\leq
\prod_{q<p}\Big(1+\frac{2}{q-1}\Big)
\prod_{\substack{q<p\\
a_q\geq 1}}
\Big(
\frac{2}{q^{a_q}}
(a_q+1)
\Big)
=
\frac{2^{\omega(g)}\tau(g)}{g}
\prod_{q<p}\Big(1+\frac{2}{q-1}\Big),
\]
where, by Mertens' estimate, the rightmost product satisfies
\[
\prod_{q<p}\Big(1+\frac{2}{q-1}\Big)\ll\log^2(2p).
\]
Substituting the bound of $\cR_p(g)$ into \eqref{def:R_p} yields \eqref{eq:kernel}, as desired.
\end{proof}

For every $g\geq 1$ with $P^+(g)<p_j$, put
\[
X(g)=\sum_{r\geq 1}
n_{\cC}(gp_j^r)p_j^{-r}.
\]
Applying Lemma~\ref{lem:kernel} to the rightmost expression in \eqref{eq:second-before-group} gives
\[
M_j^{(2)}
\leq
C_0\log^2(2p_j)
\sum_{P^{+}(g)<p_j}
X(g)^2
\frac{2^{\omega(g)}\tau(g)}{g}.
\]
By the Cauchy-Schwarz inequality,
\[
X(g)^2
\leq
\frac{1}{p_j-1}
\sum_{r\geq 1}
n_{\cC}(gp_j^{r})^2
p_j^{-r}.
\]
Since $1/(p_j-1)\leq 2/p_j$ and
$2^{\omega(g)}\tau(g)\leq \tau(g)^2 \leq \tau(gp_j^r)^2$,
we obtain
\begin{equation}
\label{eq:second-moment-stage}
 M_j^{(2)}
 \leq 2C_0
 \frac{\log^2(2p_j)}{p_j}
 \sum_{\substack{D\geq M\\P^+(D)=p_j}}
 n_{\cC}(D)^2\frac{\tau(D)^2}{D}.
\end{equation}

We now combine the two moment estimates to prove Theorem~\ref{thm:main}.

\begin{proof}[Proof of Theorem~\ref{thm:main}]
Let \[
\cP_1=\{p\in \cP:\mathsf{A}_p(m)\leq \mathsf{B}_p(m)\}
\quad
\textup{and}
\quad
\cP_2=\cP\setminus\cP_1.
\]
Recalling Remark \ref{remark1.6}, condition \eqref{eq:min-condition}
implies
\begin{equation}\label{eq:ABbounded}
\sum_{p\in \cP_1}\mathsf{A}_p(m)<\infty
\quad
\textup{and}
\quad
\sum_{p\in \cP_2}
\mathsf{B}_p(m)<\infty.
\end{equation}
For $M\geq2$, let
\begin{align*}
 S_1(M)
 &:={}
 \sum_{\substack{D\geq M\\P^+(D)\in\cP_1}}
 \frac{m(D)\tau(D)}{D},
 \\
 S_2(M)
 &:={}
 \sum_{\substack{D\geq M\\P^+(D)\in\cP_2}}
 \frac{m(D)^2\tau(D)^2
 \log^2\!\bigl(2P^+(D)\bigr)}{D\,P^+(D)}.
\end{align*}
It follows from \eqref{eq:ABbounded} that $S_1(M)$
and $S_2(M)$ are tails of convergent series with nonnegative
terms, hence
\[
S_1(M)+ 2C_0 S_2(M)\longrightarrow0
\qquad\text{as }M\longrightarrow\infty.
\]
Now, we choose $M=M(m)$ so large that
\[
 S_1(M)+ 2C_0 S_2(M)\leq\frac{1}{2}.
\]

Suppose that $\cC$ satisfies the hypotheses of Theorem~\ref{thm:main}. For $p_j\in \mathcal{P}_1$ we use \eqref{eq:first-moment-stage}, whereas for $p_j\in \mathcal{P}_2$ we use \eqref{eq:second-moment-stage}. 
Combining these bounds with $n_{\cC}(D)\leq m(D)$ and \eqref{def:M_j^12}, we obtain
\begin{equation}\label{eq:sumPB}
\sum_{j=1}^{s}
\Prob_{j}(\cB_j)
\leq
\sum_{j=1}^{s}
\min\{M_j^{(1)},
M_j^{(2)}\}
\leq S_1(M) + 
2C_0S_2(M)\leq
\frac{1}{2}.
\end{equation}
Writing
\[
R=G\setminus\bigcup_{i=1}^k C_i
 =G\setminus\bigcup_{j=1}^s\cB_j,
\]
by \eqref{invariant_prob} and \eqref{eq:sumPB}, we have
\[
\Prob_s(R)
=
1-
\Prob_s\Big(\textstyle{\bigcup_{j=1}^s\cB_j}\Big)
\geq
1-
\sum_{j=1}^{s}\Prob_j(\cB_j)
\geq
\frac{1}{2}.
\]
Finally, employing \eqref{eq:RN-bound} with $s=\omega(Q)$ implies
\[
\frac{|G\setminus \bigcup_{i}C_i|}{|G|}
\geq
 2^{-s}\Prob_s(R)
 \geq
 2^{-\omega(Q)-1}> 0,
 \]
as desired.
\end{proof}

\section{An extension to groups with a nilpotent normal subgroup}
\label{sec:bounded-index-extension}

Theorem~\ref{thm:main} extends to finite groups containing a normal
nilpotent subgroup of bounded index. The multiplicity condition
remains unchanged, while the threshold for the subgroup indices
may also depend on the bound for the index of the normal subgroup.

\begin{proposition} \label{thm:bounded-index-nilpotent}
Let $B\geq 1$ be an integer, and let $m:\{2,3,\ldots\}\longrightarrow \N\cup\{0\}$
be a function satisfying \eqref{eq:min-condition}.
Then there exists an integer $M=M(m,B)\geq 2$ with the following property.

Let $G$ be a finite group containing a normal nilpotent subgroup $N$ with $[G:N] \leq B$, and let
\[
\cC=\{C_i=a_iH_i:1\leq i\leq k\}
\]
be a finite family of left cosets of subgroups $H_i\leq G$. Set $D_i=[G:H_i]$ and $n_{\cC}(D):=
|\{i:D_i=D\}|$.
Suppose that $D_i\geq M$ for all $1\leq i\leq k$, and $n_{\cC}(D)\leq m(D)$ $(D\geq 2)$.
Then $\cC$ does not cover $G$. More precisely, define
\[
d_i=[N:N\cap H_i],
\qquad
Q_N=\lcmop(d_1,\ldots,d_k),
\qquad
Q=\lcmop(D_1,\ldots,D_k).
\]
Then $Q_N\mid Q$, and
\[
\frac{\left|G\setminus\bigcup_{i=1}^{k}C_i\right|}{|G|}
 \geq 2^{-\omega(Q_N)-1}
 \geq 2^{-\omega(Q)-1}>0.
\]
\end{proposition}

\begin{proof}
For a positive integer $B$, define
\[
m_B(d):=\sum_{e=1}^{B}m(ed),
\qquad d\geq 2,
\]
and set the constant
\[
c_B:=\sum_{e=1}^{B}e=\frac{B(B+1)}{2}.
\]
We first verify that $m_B$ satisfies the summability condition \eqref{eq:min-condition} in \Cref{thm:main}. 

Suppose first that $p>B$. By definition,
\begin{equation}\label{eq:ApmB}
\mathsf{A}_p(m_B)
=\sum_{\substack{d\geq 2\\P^+(d)=p}} \frac{m_B(d)\tau(d)}{d}
=\sum_{e=1}^{B} \sum_{\substack{d\geq 2\\P^+(d)=p}} \frac{m(ed)\tau(d)}{d}.
\end{equation}
Now, we focus on the right-hand side of \eqref{eq:ApmB}. Fixing $1\leq e\leq B$ and putting $D=ed$,
since $P^+(d)=p$ and $e<p$, we have $P^+(D)=p$. Moreover, noting $\tau(d)\leq\tau(D)$ and $1/d=e/D$, we have
\[
\sum_{\substack{d\geq 2\\P^+(d)=p}} \frac{m(ed)\tau(d)}{d}
\leq e\sum_{\substack{D\geq 2\\P^+(D)=p\\e\mid D}} \frac{m(D)\tau(D)}{D}
\leq e\,\mathsf{A}_p(m).
\]
Summing over $e$ yields $\mathsf{A}_p(m_B)\leq c_B \mathsf{A}_p(m)$. Similarly, using
\[
\left(\sum_{e=1}^{B}m(ed)\right)^2
\leq B\sum_{e=1}^{B}m(ed)^2,
\]
we obtain $\mathsf{B}_p(m_B) \leq (B\,c_B)\mathsf{B}_p(m)$.
Consequently,
\[
\min\{
\mathsf{A}_p(m_B),
\mathsf{B}_p(m_B)\}
\leq
(B\,c_B)\min\{
\mathsf{A}_p(m),
\mathsf{B}_p(m)\}
\qquad (p>B).
\]

It remains to consider the primes $p\leq B$. If
$P^+(d)=p\leq B$ and $e\leq B$, then every prime factor of $ed$
is at most $B$. The preceding change of variables therefore gives
\begin{equation}\label{eq:ApSumAq}
\mathsf{A}_p(m_B)
\leq
c_B\sum_{\substack{q\leq B\\q\ \mathrm{prime}}} \mathsf{A}_q(m).
\end{equation}
Each $\mathsf{A}_q(m)$ in the finite sum on the right-hand side
is finite. Indeed, condition \eqref{eq:min-condition} implies
that at least one of $\mathsf{A}_q(m)$ and $\mathsf{B}_q(m)$
is finite. Since $m(D)$ is a nonnegative integer, we have $m(D)\tau(D)\leq m(D)^2\tau(D)^2$. 
Consequently, if $\mathsf{B}_q(m)<\infty$, then
\[
\mathsf{A}_q(m)
\leq \frac{q}{\log^2(2q)}\mathsf{B}_q(m)<\infty.
\]
It follows from \eqref{eq:ApSumAq} that
$\mathsf{A}_p(m_B)<\infty$ for every prime $p\leq B$.
Combining the two ranges of primes, we conclude that
\[
\sum_p\min\{
\mathsf{A}_p(m_B),
\mathsf{B}_p(m_B)\}<\infty.
\]

Let $M_{\mathrm{nil}}(m_B)$ be the threshold supplied by
\Cref{thm:main} for the multiplicity function $m_B$, and set $M(m,B):=B\,M_{\mathrm{nil}}(m_B)$. Let $Y=G/N$, and let $\pi:G\longrightarrow Y$
be the quotient map. For each $i$, put $e_i=[G:NH_i]$.
Since $N\trianglelefteq G$, the set $NH_i$ is a subgroup of $G$.
The index identities
$[G:H_i]=[G:NH_i][NH_i:H_i]$ and
$[NH_i:H_i]=[N:N\cap H_i]=d_i$ imply $D_i=e_id_i$.
Since $1\leq e_i\leq[G:N]\leq B$, 
we have
\[
d_i=\frac{D_i}{e_i}
\geq \frac{D_i}{B}
\geq M_{\mathrm{nil}}(m_B).
\]

For $y\in Y$, let $F_y:=\pi^{-1}(y)$
be the corresponding fibre, and define $I_y:=\{i:C_i\cap F_y\neq\varnothing\}$.
Choose $g_y\in F_y$. Then $F_y=g_yN$, and the bijection
$n\mapsto g_yn$ identifies $N$ with $F_y$.

Suppose that $i\in I_y$, and choose $x_i\in C_i\cap F_y$. Since
$x_i\in C_i$ and $x_i\in F_y$, we have $C_i=x_iH_i$ and $F_y=x_iN$.
It follows that
\[
C_i\cap F_y
=x_iH_i\cap x_iN
=x_i(H_i\cap N).
\]
Writing $x_i=g_yn_i$ for some $n_i\in N$, we see that, under the
identification $F_y\simeq N$, the intersection $C_i\cap F_y$
corresponds to $n_i(H_i\cap N)$, and it is a left coset in $N$ of a subgroup of index $d_i$.

We next bound the multiplicity of each index in a fixed fibre. For
every $d\geq 2$,
\[
\begin{aligned}
\lvert\{i\in I_y:d_i=d\}\rvert
&\leq
\sum_{e=1}^{B}
\lvert\{i:D_i=ed,\ e_i=e\}\rvert \\
&\leq
\sum_{e=1}^{B}n_{\cC}(ed)
\leq \sum_{e=1}^{B}m(ed) =m_B(d).
\end{aligned}
\]
 Let $Q_y:=\lcmop\{d_i:i\in I_y\}$, with the convention
$Q_y=1$ if $I_y=\varnothing$.
If $I_y=\varnothing$, the whole fibre $F_y$ is uncovered.
Otherwise, Theorem~\ref{thm:main} applies to the family
$\{n_i(H_i\cap N):i\in I_y\}$ in $N$.
Since $Q_y\mid Q_N$, in both cases we obtain
\[
\frac{\left|F_y\setminus\bigcup_{i=1}^{k}C_i\right|}{|F_y|}
\geq 2^{-\omega(Q_y)-1} \geq 2^{-\omega(Q_N)-1}.
\]
Since $d_i\mid D_i$ for every $i$, we also have $Q_N\mid Q$.

Finally, all fibres $F_y$ have cardinality $|N|$. Averaging the
preceding inequality over $y\in Y$, we obtain
\[
\frac{\left|G\setminus\bigcup_{i=1}^{k}C_i\right|}{|G|}
=
\frac{1}{|Y|}
\sum_{y\in Y}
\frac{\left|F_y\setminus\bigcup_{i=1}^{k}C_i\right|}{|F_y|} \geq
2^{-\omega(Q_N)-1} \geq
2^{-\omega(Q)-1}>0.
\]
Therefore, $\cC$ does not cover $G$, which completes the proof.
\end{proof}


\begin{thebibliography}{10}

\bibitem{BBMST2022}
Paul Balister, B\'{e}la Bollob\'{a}s, Robert Morris, Julian Sahasrabudhe, and
 Marius Tiba.
\newblock On the {E}rd{\H{o}}s covering problem: the density of the uncovered
 set.
\newblock {\em Invent. Math.}, 228(1):377--414, 2022.

\bibitem{BBMST2024}
Paul Balister, B{\'e}la Bollob{\'a}s, Robert Morris, Julian Sahasrabudhe, and
 Marius Tiba.
\newblock The structure and number of {E}rd{\H{o}}s covering systems.
\newblock {\em J. Eur. Math. Soc. (JEMS)}, 26(1):75--109, 2024.

\bibitem{BNR}
Maciej Borodzik, Danny Nguyen, and Sinai Robins.
\newblock Tiling the integer lattice with translated sublattices.
\newblock {\em Mosc. J. Comb. Number Theory}, 6(4):3--26, 2016.

\bibitem{CK}
J.~E. Cremona and P.~Koymans.
\newblock Lattice coverings and homogeneous covering congruences.
\newblock {\em Trans. Lond. Math. Soc.}, to appear;
arXiv:2601.03212, 2026.

\bibitem{Erdos1950}
Paul Erd{\H{o}}s.
\newblock On integers of the form {$2^k+p$} and some related problems.
\newblock {\em Summa Brasil. Math.}, 2:113--123, 1950.

\bibitem{FFSPY2007}
Michael Filaseta, Kevin Ford, Sergei Konyagin, Carl Pomerance, and Gang Yu.
\newblock Sieving by large integers and covering systems of congruences.
\newblock {\em J. Amer. Math. Soc.}, 20(2):495--517, 2007.

\bibitem{Hough2015}
Bob Hough.
\newblock Solution of the minimum modulus problem for covering systems.
\newblock {\em Ann. of Math. (2)}, 181(1):361--382, 2015.

\bibitem{HoughNielsen2019}
Robert~D. Hough and Pace~P. Nielsen.
\newblock Covering systems with restricted divisibility.
\newblock {\em Duke Math. J.}, 168(17):3261--3295, 2019.

\bibitem{Hungerford1980}
Thomas~W. Hungerford.
\newblock {\em Algebra}, volume~73 of {\em Graduate Texts in Mathematics}.
\newblock Springer-Verlag, New York-Berlin, 1980.
\newblock Reprint of the 1974 original.

\bibitem{KKL2022}
Jonah Klein, Dimitris Koukoulopoulos, and Simon Lemieux.
\newblock On the {$j$}th smallest modulus of a covering system with distinct
 moduli.
\newblock {\em Int. J. Number Theory}, 20(2):471--479, 2024.

\bibitem{LiWangWangYi2023pre}
Huixi Li, Biao Wang, Chunlin Wang, and Shaoyun Yi.
\newblock On covering systems of polynomial rings over finite fields.
\newblock {\em Proc. Amer. Math. Soc.}, 152(9):3731--3742, 2024.

\bibitem{LiWangWangYi2025}
Huixi Li, Biao Wang, Chunlin Wang, and Shaoyun Yi.
\newblock On {E}rd{\H{o}}s covering systems in global function fields.
\newblock {\em J. Number Theory}, 266:269--280, 2025.

\bibitem{LiWangYi2023pre}
Huixi Li, Biao Wang, and Shaoyun Yi.
\newblock On the minimum modulus problem in number fields.
\newblock arXiv:2302.05946, 2023.
\newblock Revised in 2026.

\bibitem{Neumann1954}
B.~H. Neumann.
\newblock Groups covered by finitely many cosets.
\newblock {\em Publ. Math. Debrecen}, 3(3--4):227--242, 1954.

\bibitem{Sun2006}
Zhi-Wei Sun.
\newblock Finite covers of groups by cosets or subgroups.
\newblock {\em Internat. J. Math.}, 17(9):1047--1064, 2006.

\end{thebibliography}
\end{document}